\input amstex

\font\teneuf=eufm10 \font\seveneuf=eufm7 \font\fiveeuf=eufm5
\font\tenmsy=msbm10 \font\sevenmsy=msbm7 \font\fivemsy=msbm5
\font\tenmsx=msam10 \font\sixmsx=msam6 \font\fivemsx=msam5
\textfont7=\teneuf \scriptfont7=\seveneuf
\scriptscriptfont7=\fiveeuf
\textfont8=\tenmsy \scriptfont8=\sevenmsy
\scriptscriptfont8=\fivemsy
\textfont9=\tenmsx \scriptfont9=\sixmsx
\scriptscriptfont9=\fivemsx

\font\teneuf=eufm10 \font\seveneuf=eufm7 \font\fiveeuf=eufm5
\font\tenmsy=msbm10 \font\sevenmsy=msbm7 \font\fivemsy=msbm5
\font\tenmsx=msam10 \font\sixmsx=msam6 \font\fivemsx=msam5
\textfont7=\teneuf \scriptfont7=\seveneuf
\scriptscriptfont7=\fiveeuf
\textfont8=\tenmsy \scriptfont8=\sevenmsy
\scriptscriptfont8=\fivemsy
\textfont9=\tenmsx \scriptfont9=\sixmsx
\scriptscriptfont9=\fivemsx
\documentstyle{amsppt}
\NoBlackBoxes
\pagewidth{145mm}
\pageheight{220mm}

\ifx\undefined\rom
  \define\rom#1{{\rm #1}}
\fi
\ifx\undefined\curraddr
  \def\curraddr#1\endcurraddr{\address {\it Current address\/}: #1\endaddress}
\fi

\topmatter
\title A binomial coefficient identity associated to a conjecture
of Beukers
\endtitle
\rightheadtext{A binomial coefficient identity}
\author Scott Ahlgren, Shalosh B. Ekhad, Ken Ono and 
       Doron Zeilberger\endauthor
\address Department of Mathematics, Penn State University,
University Park, Pennsylvania 16802 \endaddress
\email ahlgren\@math.psu.edu \endemail

\address Department of Mathematics, Temple University,
Philadelphia, Pennsylvania 19122\endaddress
\email ekhad\@math.temple.edu; http://www.math.temple.edu/
 \~{}ekhad \endemail
 
\address Department of Mathematics, Penn State University,
University Park, Pennsylvania 16802 \endaddress
\email ono\@math.psu.edu; http://www.math.psu.edu/ono/  \endemail
 
\address Department of Mathematics, Temple University,
Philadelphia, Pennsylvania 19122
\endaddress
\email zeilberg\@math.temple.edu; http://www.math.temple.edu/\~{} zeilberg
\endemail
\thanks 
The  third author is supported by NSF grant DMS-9508976 and
NSA grant MSPR-Y012. The last author is supported in part by
the NSF.
\endthanks
 
\endtopmatter

\document
If $n$ is a positive integer, then let
  $\dsize{A(n):=\sum_{k=0}^{n}{n\choose k}^2{n+k\choose k}^2}$,
and define integers $a(n)$ by
  $\dsize{\sum_{n=1}^{\infty}a(n)q^n:=q\prod_{n=1}^{\infty}
   (1-q^{2n})^4(1-q^{4n})^4=q-4q^3-2q^5+24q^7- \cdots.}$ 
Beukers conjectured that
if $p$ is an odd prime, then
$$
  A\left ( \frac{p-1}{2}\right ) \equiv a(p) \pmod{p^2}. \tag{1}
$$

In [A-O] it is shown that (1) is implied by the truth
of the following identity.
\proclaim{Theorem} If $n$ is a positive integer, then
$$
  \sum_{k=1}^{n}k{n\choose k}^2{n+k \choose k}^2
  \left \{ \frac{1}{2k}+\sum_{i=1}^{n+k}\frac{1}{i}+
 \sum_{i=1}^{n-k}\frac{1}{i}-2\sum_{i=1}^{k}\frac{1}{i}\right \}=0.
$$
\endproclaim
 
\medskip
 
\noindent
{\bf Remark.} 
This identity is easily verified using the WZ method, in a generalized form
[Z] that applies when the summand is a hypergeometric term times a WZ
potential function.
It holds for all positive $n$,
since it holds for $n$=1,2,3 (check!), and since the
sequence defined by the sum satisfies a certain 
(homog.) third order linear recurrence
equation. To find the recurrence, and its proof,
download the Maple package EKHAD and the Maple program zeilWZP
from {\eighttt http://www.math.temple.edu/\~{} zeilberg }.
Calling the quantity inside the braces $c(n,k)$, compute
the WZ pair $(F,G)$, where $F=c(n,k+1)-c(n,k)$ and
$G=c(n+1,k)-c(n,k)$. 
Go into Maple, and type 
{\eighttt read zeilWZP;
zeilWZP(k*(n+k)!**2/k!**4/(n-k)!**2,F,G,k,n,N):
}

\enddocument